# The Minimum Latency Problem

Avrim Blum [*]    Prasad Chalasani [†]    Don Coppersmith [‡]
Bill Pulleyblank [‡]    Prabhakar Raghavan [‡]    Madhu Sudan [‡]

## Abstract

We are given a set of points $p_1, \ldots, p_n$ and a symmetric distance matrix $(d_{ij})$ giving the distance between $p_i$ and $p_j$. We wish to construct a tour that minimizes $\sum_{i=1}^{n} \ell(i)$, where $\ell(i)$ is the *latency* of $p_i$, defined to be the distance traveled before first visiting $p_i$. This problem is also known in the literature as the *deliveryman problem* or the *traveling repairman problem*. It arises in a number of applications including disk-head scheduling, and turns out to be surprisingly different from the traveling salesman problem in character. We give exact and approximate solutions to a number of cases, including a constant-factor approximation algorithm whenever the distance matrix satisfies the triangle inequality.

[*]School of Computer Science, CMU. Supported in part by NSF National Young Investigator grant CCR-9357793.

[†]School of Computer Science, CMU.

[‡]IBM T.J. Watson Research Center.

## 1. Introduction

Consider a server (a repairman or a disk head, perhaps) facing a set of requests, each of which is a point. The server must schedule its visits so as to minimize the average time the requests wait before being visited. (We assume that a request is serviced instantaneously when visited.) This is a simple and natural combinatorial optimization problem faced in many day-to-day situations, and may be formalized as follows.

> Given a set of $n$ points, a symmetric distance matrix $(d_{ij})$, and a tour which visits the points in some order, let the *latency* of a *point* $p$ be the length of the tour from the starting point to $p$. Let the *total latency* of the *tour* be the sum of the latencies of all the points. We wish to find the tour which minimizes the total latency. *Of particular interest is the case when the distance matrix satisfies the triangle inequality.*

We henceforth abbreviate the minimum latency tour by MLT. Whenever we speak of approximating the MLT, we seek to find a tour whose cost approximates that of the MLT.

At first glance the MLT problem seems to be a simple variant of the traveling salesperson problem (TSP). However, closer examination reveals a variety of aspects in which this problem is very different from the TSP. Small changes in

0

the structure of a metric space can cause highly non-local changes in the structure of the MLT (see e.g., Fig 1). As a result, a simple paradigm which applies to most optimization problems on graphs — namely, " decompose the graph into biconnected components; solve the optimization problem on the various (hopefully) small components and then put them together to obtain a solution for the big problem" — does not work here. The absence of this paradigm rules out simple algorithms even in the case where the underlying graph is a tree (i.e, the metric space is simply the metric closure of a weighted tree.) Another prominent difference between the MLT and the optimal TSP tour is that the MLT may revisit points an unbounded number of times even when the underlying graph has a bounded degree. Consider the case where the metric space is the real line and the point $p_i$ is located at $(-3)^i$. In this case the MLT starting at the origin will visit the points in order (i.e., the $i$th point to be visited is $p_i$). Thus the MLT crosses the origin $n-1$ times! By a simple perturbation to this set of points, we get an example of points on the plane where the MLT is not planar (!) — again a phenomenon which distinguishes it from an optimal TSP tour.

However, it is easy to see that this problem is at least as hard as the TSP. For instance, given a set of points on which we wish to minimize the length of the TSP tour, one can augment the set of points with $N$ points at "infinity" (for a large number $N$), so that the MLT on the augmented set of points will have to minimize the length of the TSP on the original set of points. This connection shows that the MLT problem is NP-hard even in the case where the metric space is a plane.

For general metric spaces, it is possible to reduce the TSP where all distances are either 1 or 2 (a metric space) to the MLT. In conjunction with the MAX SNP-hardness of the TSP(1,2) problem [12] and the non-existence of polynomial time approximation schemes for MAX SNP-hard problems [2] this connection implies that one cannot hope for arbitrarily good approximations to the minimum latency problem on general metric spaces. Lastly, as in the case of the TSP, the MLT is NP-hard to approximate to within any bounded ratio when the distance function is arbitrary (i.e., does not satisfy the triangle inequality.)

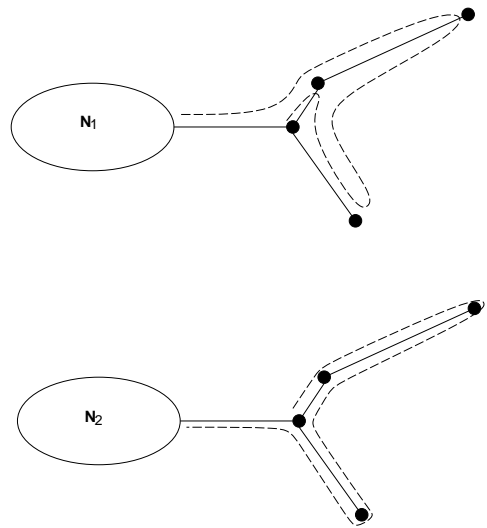

Figure 1: Non-local change in tour, due to $N_1 << N_2$

Both the MLT problem and the TSP problems are special cases of a more general problem, the "time-dependent traveling salesman problem". Here the distance function on the metric space varies with "time", i.e., if the $i$th edge to be traversed has weight $e$ then it incurs a cost $c(e, i)$. The goal is to minimize the total cost of visiting all vertices. The TSP is the special case where the cost function is just $c(e, i) = e$ (independent of $i$) and the MLT is the case where the cost function is given by $c(e, i) = (n - i) * e$. Situations where the cost function $c(e, i)$ is some convex combination of these two seem to arise naturally too. For instance, consider the following vehicle routing problem: a delivery truck has to deliver $N$ items to $N$ points in a metric space, and then return to its starting point. If the truck travels a distance $D$ with $k$ items, then the cost of that leg is $(k + u)D$ (the cost is proportional to the load on the truck which is $k$ units for the items $+ u$ units for the weight of the truck).

The time-dependent TSP and the MLT problem have been studied earlier, under the names of *the deliveryman problem* or *the traveling repairman problem* [10, 9, 13, 5]. In [10] and [13] it is shown that any depth-first route is an optimal MLT on a tree with unit edge weights. Minieka also gives an exponential time algorithm for solving the problem optimally on general trees. In [9] and [5] exponential time algorithms are provided to compute the MLT optimally in general metric spaces. Afrati *et al.* [1] study a version of the problem in which some requests have prescribed deadlines on their latencies. Among other results, they show that the latency problem on the line can be solved in polynomial time by dynamic programming.

## 1.1. Our Results

Given the evidence against the existence of polynomial time algorithms to solve the exact problem in general cases, we first isolate some situations in which the MLT problem can be solved exactly in polynomial time. We start by proving that any depth first search tree gives an optimal MLT for the unweighted tree (all edge lengths are 1). This theorem has been proved in the earlier mentioned works of [10, 13], but our proof is simpler. Surprisingly, this theorem does not hold for weighted trees (i.e., depth first search trees are not optimal traversals). We consider a special case of weighted trees — diameter 3 trees — and give a simple dynamic programming algorithm to solve this case optimally.

For general metric spaces, we provide a constant factor approximation algorithm in Section 3. This approximation algorithm is based on the following inequality which we believe may be of independent interest.

Let $p$-latency denote the minimum total latency over all tours which start at $p$. Let the $i$-tree$(p)$ denote the weight of the shortest tree spanning $i$ nodes including $p$. Then

$$\sum_{i=1}^{n} i\text{-tree}(p) \leq p\text{-latency} \leq 8 * \sum_{i=1}^{n} i\text{-tree}(p)$$

The above inequality is constructive in the sense that given an algorithm which finds a tree on $i$ nodes with weight approximately $i$-tree$(p)$, we can obtain an algorithm to obtain a tour starting at $p$ which approximates the MLT. Unfortunately, the task of finding the best $i$-tree is also NP hard [6] for general metric spaces and hence the above inequality does not solve the problem for us. We are able to find an exact algorithm for the $i$-tree problem for the case of the weighted tree. Garg and Hochbaum [6] also give an $O(\log i)$ approximation algorithm for the case when the metric space is Euclidean in a fixed number of dimensions, and this yields a $O(\log n)$ approximation to the MLT in such spaces.

To obtain a constant factor approximation algorithm for general metric spaces, we show that an inequality similar to the above holds even when the algorithm provides only a solution with $i'$ vertices with cost within a constant factor of the best $i$-tree. (The exact relationship between $i$ and $i'$ is described in Section 3.2.) In fact, we show that it suffices to do this for values of $i$ close to $n$. We then use a constant-factor approximation algorithm to the prize-collecting TSP [8] to obtain an algorithm extending the above ideas.

The solution we find to approximate the MLT also turns out to approximate the TSP simultaneously. This allows us to extend our results on the MLT approximation to include a family of time-dependent TSP problems which includes the truck delivery problem described earlier. The family of time dependent TSPs we consider is the case when the cost function on an edge is linear in $i$. This is described in Section 3.3.

## 2. Exact Solutions

In this section we give polynomial-time algorithms for computing the MLT in some special cases.

### 2.1. Unweighted trees and depth-first search

Consider the case when the points are vertices of a tree all of whose edges have unit length. Minieka [10] and Reynolds [13] have shown that a tour is optimal if and only if it is a depth-first search. Here we give a very simple and succinct proof of this result; we proceed by showing that every depth-first search is optimal, and that optimal tours are depth-first searches.

For a vertex $v$, let $depth(v)$ denote its distance in the tree from the starting point. The following two claims imply the result.

**Claim 1:** For any tour, the $i$th distinct vertex to be visited has latency at least $2i - depth(v)$.

**Claim 2:** Fix a depth-first search, and let $v$ be the $i$th distinct vertex to be visited. Then $v$ has latency exactly $2i - depth(v)$.

**Proof:** Fix some index $i$. Both claims are easy to see by imagining that the walk returns to the origin after visiting the $i$th distinct vertex. It is clear that after doing so, any such walk must have traveled at least $2i$ steps. Furthermore the DFS achieves this exactly, since a DFS never traverses any edge more than twice. But when the walk is forced to return to the origin, this extends the length of the walk by exactly $depth(v)$. Thus any walk has a latency at $v$ of at least $2i - depth(v)$ and the DFS achieves this bound exactly.   □

### 2.2. Dynamic programming

Dynamic programming may be used to obtain optimal solutions in polynomial time for cases when there is a good bound on the number of potential partial solutions. We illustrate this in two cases: when the $p_i$ are points on the line, and when they are vertices of a tree of diameter 3. Although Afrati et al. [1] have already shown a dynamic programming algorithm for the line, we include here a brief outline for completeness.

When the $p_i$ are points on the line, a partially complete tour covers an interval of the line that includes the starting point. A *state* in the dynamic programming tableau is an interval of the line, together with one end of the interval; it represents our current position at the end of a partial tour (it suffices to consider intervals whose end-points are both points $p_i$ from the input). Clearly there are $O(n^2)$ states. The algorithm works from large intervals down to smaller ones. For each state it takes constant time to decide whether the first point visited by the MLT starting in this state is to the left or to the right of the covered interval. Thus we have:

**Theorem 1 ([1]):** *When the $p_i$ are points on the real line, dynamic programming yields an optimal solution in $O(n^2)$ time.*

Dynamic programming also yields an optimal solution in $O(n^2)$ time when the points are vertices of a tree of diameter 3, when there are positive real lengths on the edges. Such a tree consists of a central edge with spokes hanging off its two end-points. The observation here is that the spokes at either end are visited in increasing order of their lengths. Thus a state for dynamic programming consists of two integers $k_L$ and $k_R$, where $k_L$ is the number of spokes visited at the left end-point of the central spoke, and $k_R$ the number of spokes at the right end-point (in addition, we also use a bit to keep track of whether we're at the left or the right hub). Again, the number of states is $O(n^2)$. Unfortunately, this technique does not seem to extend even to trees of diameter 4.

## 3. Approximation Algorithms

In this section we focus on techniques for devising approximation algorithms. We begin with a familiar algorithm on the real line: starting at the origin, we are to reach a point $p$ at an unknown position on the line. Let $d$ be the distance from the origin to $p$. Baeza-Yates, Culberson and Rawlins [3] have given a simple deterministic "doubling" algorithm that finds $p$ after walking at most a distance $9d$. By adapting this algorithm, we have a simple deterministic algorithm that achieves a 9-approximation for points on the line (the dynamic programming approach above gives the exact solution, but is more computationally intensive). One can extend this approach to the case when the $p_i$ are vertices of a *layered graph*, adapting familiar algorithms for layered graph traversal [11] to get provably good approximations. We omit these details here.

### 3.1. Reduction to $i$-trees

Given $n$ points and a distance metric between them, the *$i$-tree problem* is the following: we are to find the shortest tree spanning $i$ of the input points. We now show that approximating the MLT reduces to solving (or approximating) the $i$-tree problem on the same input points, for all $1 \leq i \leq n$. The reduction works provided the distances satisfy the triangle inequality. The $i$-tree problem is NP-hard to solve in general, but can be $O(\log n)$-approximated in Euclidean spaces of fixed dimension [6].

To simplify the presentation, let us assume that the distance from the starting point to its nearest neighbor is 1. Renumber the points 1 through $n$, where the point $i$ is the $i$th vertex to be visited by the optimal tour. Let $S_j$ be the set of vertices with latency between $2^j$ and $2^{j+1}$ in the optimal tour, and let $n_j = |S_j|$.

By invoking an algorithm for the $i$-tree problem for $i = 1, 2, \ldots$, we can determine, for each $j$, the maximum number of points that we can span by a tree of cost at most $2^{j+1}$; denote this number by $m_j$. Our approximation to the MLT is now the following: for $j = 1, 2, \ldots$, we traverse the $m_j$-tree (say in depth-first fashion), returning to the origin between successive values of $j$. Consider the $i$th vertex visited by such a tour, and let $i \in S_j$. Then the latency of the $i$th vertex in the optimal tour is at least $2^j$. We now show that the latency of the $i$th vertex in our tour is at most $8 * 2^j$. Since there exists a tree of size $2^{j+1}$ which visits at least $\sum_{k \leq j} n_j \geq i$ points, $m_j$ is at least $i$. Thus the latency of the $i$th vertex in our tour is at most $2 * (\sum_{k < j} 2^{k+1} + 2^{j+1}) \leq 8 * 2^j$. Thus the total latency of our tour is at most 8 times the latency of the optimal MLT.

It is easy to see that the above argument extends even if we only have a $c$-approximation algorithm to the $i$-tree problem. Thus we get the following:

**Theorem 2:** *Given an algorithm that gives a $c$-approximate solution to the $i$-tree problem, we can obtain a $8c$-approximation algorithm to the minimum latency problem.*

We conclude with an example of a case where the $i$-tree problem can be solved exactly. Suppose that the input points are vertices of a tree whose edges have positive real lengths. In this case, dynamic programming can be used to obtain $i$-trees of optimal length as follows. First, suppose the graph were a binary tree. In this case, for each point, the optimal solution to the $i$-tree problem on the subtree rooted at that point can be computed given solutions to the $i'$-tree problem for the point's children for all $i' \leq i$. This leads immediately to a dynamic programming solution for *binary* trees. Notice that this solution can be extended easily to the more general problem in which points have weights in $\{0, 1\}$ and the goal is to find the shortest tree with a total weight of points $i$. Now, to handle non-binary trees, simply replace a vertex of higher degree with a binary tree whose root has weight 1, all other vertices have weight 0, and the edges have weight 0 as well. Thus our reduction implies an 8-approximation to the minimum latency problem in this case.

## 3.2. A constant-factor approximation for metric spaces

We now describe a constant factor approximation algorithm for the minimum latency problem in general metric spaces. A key point here is that it is sufficient to approximate the $i$-trees problem only when $i$ is quite large: in order to achieve a constant-factor we need only concern ourselves with the latency of the last fraction of vertices visited. This insight motivates the following definition:

**Definition:** An $(\alpha, \beta)$-TSP-approximator is an algorithm that given bounds $\epsilon$, $L$, an $n$-point metric space $M$, and a starting point $p$ finds a tour starting at $p$ of length at most $\beta L$, which visits at least $(1 - \alpha\epsilon)n$ vertices, if there exists a tour of length $L$ which visits $(1 - \epsilon)n$ vertices.

**Comment:** If only one of the two parameters among $\epsilon$ and $L$ is specified, we can perform a binary search for the optimal value of the other parameter. Hence in what follows we will sometimes call the approximator with only one of the parameters specified.

We first show how to construct a $(3, 6)$-TSP-approximator from a 2-approximation algorithm for the "prize-collecting traveling salesman problem" due to Goemans and Williamson [8] (which improves on a 5/2-approximation due to Bienstock et al. [4]). We then describe how the TSP-approximator is used to approximate the latency problem.

The prize-collecting traveling salesman problem is the following.

> Given a weighted graph with penalties on the vertices, the cost of a tour on some subset of the vertices is the total distance traveled plus the sum of the penalties on the vertices *not* visited. We wish to find the tour, beginning at some prespecified root vertex, that minimizes this cost.

**Lemma 3 ([8]):** *There exists a $2 - \frac{1}{n-1}$ approximation algorithm for the prize collecting traveling salesman problem, on any graph which satisfies the triangle inequality.*

From this we obtain the following implication.

**Corollary 4:** *There exists a $(3, 6)$-TSP-approximator.*

**Proof:** Let $\epsilon$ be the fraction of vertices missed out by the optimal tour of length $L$. Let us begin by assuming $\epsilon$ is known. Place a penalty $P = 2L/(\epsilon n)$ on each vertex. So, there exists a tour of total penalty at most $L + \epsilon n P = 3L$ in the prize-collecting sense. Therefore, the GW algorithm finds a tour starting at the point $p$, that visits at least a $(1 - 3\epsilon)$ fraction of the vertices (else the penalty is more than $3\epsilon n P = 6L$) and of total distance $6L$. If $\epsilon$ is not known then just perform binary search, selecting the tour found of length less than $6L$ that visits the most vertices. □

We now show how the $(\alpha, \beta)$-TSP-approximator can be used to find a tour with small latency.

**Lemma 5:** *A tour of latency at most $8\lceil\alpha\rceil\beta$ times the optimal can be found by making polynomially many calls to an $(\alpha, \beta)$-TSP-approximator.*

**Proof:** The proof of this lemma mimics the $i$-tree reduction closely. For simplicity of presentation, we assume that the distance from the starting point to its nearest neighbor is 1. We also assume $\alpha$ is integral (or else we can use its ceiling as our bound).

Our approximation algorithm calls the approximator for $L = 2, 4, 8, \ldots, 2^i, \ldots$ and then concatenates these tours to obtain a tour which visits all the vertices.

For the analysis of this tour, we partition the vertices of our tour into blocks of size $\alpha$, where the $i$th block $B_i$ contains the vertices which are the $n - \alpha(i + 1) + 1, n - \alpha(i + 1) + 2, \ldots, (n - \alpha i)$th to be visited by our tour. Now consider the minimum latency tour, and let $S_j$ be the set of vertices of latency between $2^{j-1}$ and $2^j$. Let the $n - i$th vertex visited by the optimal tour lie in the set $S_j$ and let $i = (1 - \epsilon)n$. Consider the $\alpha$ vertices in the block $B_i$. All these vertices must have been visited by our $j$th round trip, and thus the latency in our tour for each of these vertices is bounded by $\beta(2 * 2^j + 2 * 2^{j-1} + \cdots) \leq \beta 2^{j+2}$. Thus the latency of these $\alpha$ vertices in our tour is at most $\alpha\beta 2^{j+2}$, which is at most $8\alpha\beta$ times the latency of the $i$th vertex. We can thus charge the cost of each block of vertices in our tour against the latency of distinct vertices in the optimal MLT, with a multiplicative constant of at most $8\alpha\beta$.

□

**Theorem 6:** *There is a polynomial-time 144-approximation algorithm for the MLT whenever the distances $(d_{ij})$ satisfy the triangle inequality.*

**Proof:** Follows from Lemma 5 and Corollary 4. □

The constant of 144 above can be improved easily to 72 using the following idea: We call the $(\alpha, \beta)$-approximator for $\epsilon = 1/2, 1/4, \ldots, 2^{-i}, \ldots,$. We concatenate the tours so obtained (in the order they are constructed). An analysis similar to that in Lemma 5 shows that this gives a $4\alpha\beta$ approximation to the MLT.

**Theorem 6':** *There exists a polynomial time 72-approximation algorithm for the MLT whenever the distances $(d_{ij})$ satisfy the triangle inequality.*

This constant of 72 has been further improved by recent work of Goemans and Kleinberg [7] who give a factor of 29 approximation algorithm for the MLT. However the 144 approximation algorithm given by Theorem 6 has some additional properties which we use in the next section.

### 3.3. Simultaneous Approximation of the TSP and MLT

Consider the length of the tour produced by our MLT approximation in Theorem 6 above. Could it be the case that this length is much longer than the optimal TSP tour? We give a negative answer to this question. In the process we extend our technique to approximate a class of problems that combine the objectives of the MLT and TSP in any convex combination: namely the positive-linear time-dependent TSP as defined next.

> Given a set of $n$ points, a symmetric distance matrix $(d_{ij})$ and a time-dependent cost function $c : E \times \mathcal{Z}_+ \to \mathcal{R}$, the *time dependent traveling salesman problem (TDTSP)* is that of finding a tour which visits each vertex exactly once and minimizes $\sum_{i=1}^{n} c(e_i, i)$, where $e_i$ is the $i$th edge used by the tour. When the cost function is of the form $c(e, i) = (a * i + b) * d_e$, then this is the *linear TDTSP*. If $c(e, i) = (a * (n - i) + b) * d_e$ for positive $a, b$ and $d_e$'s, then this is the *positive-linear TDTSP*.

We start by observing that the above algorithm not only produces a tour which approximates the MLT but also the TSP length. This is seen as follows: Let the length of the longest tour above be $L$. Then $L$ is at most $2\beta$ times the length of the TSP. Furthermore the length of our tour is at most $\sum_{j \geq 0} \frac{L}{2^j}$, which is at most $2L$. Thus the tour above is actually a $4\beta$ approximator to the minimum TSP as well. This allows us to obtain the following approximation for the class of positive-linear TDTSP problems.

**Theorem 7:** *There is a polynomial-time 144-approximation algorithm for any positive-linear TDTSP.*

**Proof:** Consider a cost function of the form $c(e,i) = d_e * (a * (n - i) + b)$. The cost of the optimum tour for this cost function is at least $a * \text{opt(MLT)} + b * \text{opt(TSP)}$ (where opt(MLT) and opt(TSP) represent the optimal lengths of the MLT and TSP on the graph with distance function $d_e$). We use the algorithm described above to produce a tour of length at most $24 * \text{opt(TSP)}$ and latency at most $144 * \text{opt(MLT)}$. Thus our cost on this cost function is at most $144 * (a * \text{opt(MLT)} + b * \text{opt(TSP)})$, which is within a factor of 144 of the optimal. □

It may be pointed out that the above theorem also gives a 144 approximation algorithm for the case when the cost function is of the form $c(e,i) = d_e * (a*i+b)$ for positive $a$ and $b$ (we simply reverse the tour returned by the above procedure). Another case of linear TDTSP for which a constant factor approximation algorithms exists is the case where both $a$ and $b$ are negative (negative-linear TDTSP). In this case minimizing this negative function is the same as maximizing the TSP and the latency of the tour. It turns out the "greedy tour", i.e., one that visits the furthest unvisited vertex next, gives a factor of 2 approximation to both the MaxLT and the TSP, thus giving a factor of 2 approximation to the negative-linear TDTSP. Details of this proof will be provided in the full paper.

## 4. Conclusion

Very recently, Goemans and Kleinberg [7] have shown that our factor of 72 approximation for the MLT can be improved to a factor slightly less than 29. It remains to be seen whether can be used to improve the constant for the positive-linear TDTSP. The question of whether there exists an exact solution to the MLT for the case of weighted trees remains open as well.

The online version of the latency seems to be a very interesting case for further study. Here the problem would be formalized as follows: "The input is a sequence of requests which arrive online: A request arriving at time $t$ and getting serviced at time $t'$ incurs a latency of $t' - t$. We wish to minimize the average latency of all requests." Once again a variety of applications can be found for this problem e.g. diskhead scheduling, parcel pickup etc. The online considerations introduce a variety of forms of algorithms that we may wish to consider and we discuss some of these models here.

We start by observing that no competitive algorithms can be found for this problem when competitivity is measured against the best offline algorithm for a specific request sequence. For instance, consider the case of a two-point metric space where the points are separated by a unit distance and requests are generated randomly at one of the two points every unit of time. The offline algorithm does not incur any cost while any online algorithm incurs a latency of 1 with probability 1/2 every step, giving us the following consequence:

**Proposition 8:** *There is no on-line latency algorithm achieving any bounded competitiveness, even against an oblivious adversary.*

A more reasonable comparison is that of our algorithm against other online algorithms in the model where the requests are generated independently at every time step according to some distribution $\mathcal{D}$ on the points. For the special case of the complete graph, (i.e., the distance between every pair of points is 1), the following can be shown:

**Proposition 9:** *For the complete graph on $n$ points, there exists a simple online algorithm which achieves a competitive factor of 2 against any online algorithm when the requests are drawn independently at random according to any distribution $\mathcal{D}$.*

Of course, the case of general metric spaces is completely open.